\documentclass{article}
\usepackage{amssymb}
\usepackage{romp33}
\addtocounter{section}{-1}
\newtheorem{defi}{Definition}
\newtheorem{nota}{Notations}
\newtheorem{lem}{Lemma}
\newtheorem{thm}{Theorem}
\newtheorem{rem}{Remark}

\newcommand{\lbr}{\linebreak[0]}
\renewcommand{\phi}{\varphi}

\def\G{\mathit{\Gamma}}

\def\g{\gamma}
\def\ra{\mathop{\rm rank}\nolimits}

\def\Re{\mathop{\rm Re}\nolimits}
\def\Im{\mathop{\rm Im}\nolimits}
\def\mathclass#1{\setbox\mathbox=\hbox{\vfootnote{}}}
\def\Id{\mathop{\rm Id}\nolimits}

\def\R{{\Bbb R}}
\def\C{{\Bbb C}}

\def\P{{\Bbb P}}
\renewcommand{\le}{\leqslant}
\renewcommand{\ge}{\geqslant}

\begin{document}
\title{On integrability of generalized Veronese curves of
distributions}
\author{Andriy
Panasyuk\thanks{Partially
supported by the Polish grant KBN 2 PO3A 135 16.}
\\
Division of Mathematical Methods in
Physics,\\ University of Warsaw,\\ Ho\.{z}a St.~74, 00-682
Warsaw, Poland,\\ e-mail:  panas@fuw.edu.pl\\ and \\
Pidstrygach Institute for Applied Problems of \\
Mathematics
and Mechanics,\\ Naukova Str. 3b,\\ 79601 Lviv, Ukraine}
\maketitle
\begin{abstract}
Given a 1-parameter family of 1-forms  $\g(t)=
\g_0+t\g_1+\ldots+t^n\g_n$, consider
the condition $d\g(t)\wedge\g(t)=0$ (of
integrability for the annihilated by $\g(t)$ distribution
$w(t)$).
We prove that in order that this condition
is satisfied for any $t$ it is sufficient
that it is satisfied for $N=n+3$ different values of $t$
(the corresponding implication for $N=2n+1$ is obvious).  In
fact we give a stronger result dealing with distributions
of higher codimension.  This result is related to the
so-called Veronese webs and can be applied in the theory of
bihamiltonian structures.
\end{abstract} \date{}

\section{Introduction}

The notion of a Veronese web was introduced by I.M.Gelfand
and I.S.Zakharevich \cite{gz2} as a natural
invariant of bihamiltonian structures of corank 1. They
conjectured that locally this invariant is complete, i.e.
determines the bihamiltonian structure up to an
isomorphism. This conjecture was proved by F.J.Turiel
\cite{t1}.

Let us briefly recall relevant definitions. Assume that we
have a 1-parameter family $\{w(t)\}_{t\in \R\P^1},
w(t)\subset TM,$ of distributions of codimension 1 on a
smooth manifold $M$ such that in a neighbourhood of any
point there exist an annihilating $w(t)$ 1-form
$\g(t)\in\G(w(t))^\bot$ and a coframe
$\{\g_0,\ldots,\g_n\}$ with the property
$\g(t)=\g_0+t\g_1+\cdots+t^n\g_n$ (we assume that
$\g(\infty)=\g_n$ annihilates $w(\infty)$ and write $\G$
for the space of sections of a vector bundle). We call
this family of distributions a Veronese curve of
distributions. We say that it is integrable or it is a
Veronese web if each distribution $w(t)$ is integrable,
i.e.
\begin{equation} d\g(t)\wedge\g(t)=0, t\in \R
\label{1}
\end{equation}
and
\begin{equation}
d\g(\infty)\wedge\g(\infty)=0
\label{2}
\end{equation}
(it can be shown that $(\ref{1})\Rightarrow(\ref{2})$).

The equation (\ref{1}) is polynomial of degree $2n$ in $t$,
hence it is sufficient that it is satisfied for $2n+1$
different values of $t$ in order that it is satisfied for
any $t\in \R$. In other words, integrability of $w(t)$ at
$2n+1$ different points implies integrability of a Veronese
curve of distributions.

It is remarkable that this number can be essentially
reduced. I.Zakharevich \cite{z2} conjectured
that the integrability of $w(t)$ at $n+3$ different
points implies the integrability of the Veronese curve of
distributions.

The aim of this paper is to prove this result in a more
general setting. We give it for generalized Veronese curves
of distributions.  The precise definition is given in
Section 1. Here we shall only mention that integrable
generalized Veronese curves of distributions coincide with
a particular case of the so called Kronecker webs
introduced in \cite{z1}. The last notion serves as an
invariant of bihamiltonian structures of higher corank in
the same manner as the notion of a Veronese web does this
for bihamiltonian structures of corank 1 (see also
\cite{t2}, \cite{p}).

One of possible applications of our result is that it
allows to generalize theory of nonlinear wave equations
~\cite{z2} to higher dimensional and codimensional cases.
Also it may help to establish new relations between
Veronese webs and classical webs. In ~\cite{n} such
relations are studied in the case $n=1$.

Briefly, our method can be described as follows. Given
$n+3$ foliations corresponding to the integrable points of
the curve of distributuins $w(t)$, we construct locally a
curve of foliations $\widetilde{w}(t)$ in a larger space
such that $\widetilde{w}(t)$ projects onto $w(t)$, so
proving the integrability of $w(t)$.  Our proof works only
in the real-analytic category, since it uses the
complexifications of the objects. This, of course, is one
of the disadvantages of the method.

The author wishes to thank Ilya Zakharevich for helpful
discussions.

All objects in this paper are from $C^\omega$-category.

\section{Basic definitions}

\begin{defi}
\label{d1}
Let $M^{k(n+1)}$ be a manifold of dimension $m=k(n+1)$ and
let $\{ w(t)\}_{t\in\R\P^1},\linebreak w(t)\subset TM$, be
a family of distributions of codimension $k$ (as subbundles
they have rank $kn$). Assume that in a neighbourhood of any
point there exist $k$ independent annihilating $w(t)$
1-forms $\g^1(t),\ldots,\g^k(t)\in \G(w(t))^\bot$ and a
coframe
$\{\g_0^1,\ldots,\g_n^1,\g_0^2,\ldots,\g_n^2,\ldots,\lbr
\g_0^k,\ldots,\g_n^k\}$ such that
$\g^i(t)=\g_0^i+t\g_1^i+\cdots+t^n\g_n^i,i\in\overline{1,k}$. Then we call
$\{ w(t)\}_{t\in\R\P^1}$ a generalized Veronese curve of
distributions.
\end{defi}
\begin{defi}\label{d2}
A generalized Veronese curve of distributions $\{
w(t)\}_{t\in\R\P^1}$ is called integrable or a generalized
Veronese web if each distribution $w(t)$ is integrable,
i.e.
\begin{equation}
\label{e3}
d\g^i(t)=\sum_s\beta^i_s(t)\wedge\g^s(t), i\in\overline{1,k}
\end{equation}
for some depending on t 1-forms $\beta_s^i(t),
i,s\in\overline{1,k}$ (we put
$\g^i(\infty):=\g_n^i,i\in\overline{1,k}$). The
corresponding foliations will be denoted by ${\cal W}(t)$.
\end{defi}
\begin{rem}\rm It can be proved that if equations
(\ref{e3}) are  satisfied for each $t\in\R$, then
$w(\infty)$ is automatically integrable.
\end{rem}
\begin{rem}\rm
In case $k=1$ one gets the standard definition of a
Veronese web. For $k>1$ we obtain a particular case of
Kronecker webs (\cite{z1}) with Kronecker blocks of equal
dimension. In \cite{p} the terminology "generalized
Veronese webs" was used for slightly different
objects, namely for Kronecker webs without Kronecker blocks
of equal dimension.
\end{rem}
\begin{rem}
\rm
If $\{
w(t)\}_{t\in\R\P^1}$ is a generalized Veronese web and
$a_1,\ldots, a_l\in \R\P^1, l\ge m$, are different, then
$\{{\cal W}(a_1),\ldots{\cal W}(a_l)\}$ is a "classical"
$l$-web of codimension $k$ (see \cite{g}), i.e. the
foliations ${\cal W}(a_1),\ldots{\cal W}(a_l)$ are in
general position.
\end{rem}

\section{Main theorem}

\begin{lem}
Let $M$ be a complex manifold, $J:TM\rightarrow TM$ be the
complex structure operator on the real tangent bundle $TM$.
Assume that $F\subset TM$ is an integrable distribution
such that the distribution $JF\subset TM$ is also
integrable. Then the distribution $(\Id+tJ)F$ is integrable for
any $t\in\R$.
\end{lem}

\noindent{\sc Proof} Let $\{v_1,\ldots,v_s\}$ be a
local system of generating $F$ vector fields. Then
$$[v_i,v_j]=\sum_l\alpha_{ij}^lv_l,
[Jv_i,Jv_j]=\sum_l\beta_{ij}^lJv_l$$
for some functions $\alpha_{ij}^l,\beta_{ij}^l$. The
following calculations use the integrability condition for
the complex structure
$$J[v_i,v_j]-J[Jv_i,Jv_j]=[Jv_i,v_j]+[v_i,Jv_j]$$
and complete the proof:
$$
[v_i+tJv_i, v_j+tJv_j]=$$
$$[v_i,v_j]+t^2[Jv_i,Jv_j]+t([v_i,Jv_j]+[Jv_i,v_j])= $$
$$[v_i,v_j]+t^2[Jv_i,Jv_j]+t(J[v_i,v_j]-J[Jv_i,Jv_j])=$$
$$[v_i,v_j]+tJ[v_i,v_j]+t(-J[Jv_i,Jv_j]+t[Jv_i,Jv_j])= $$
$$(\Id+tJ)[v_i,v_j]+t(\Id+tJ)(-J[Jv_i,Jv_j])=$$
$$\sum_l\alpha_{ij}^l(v_l+tJv_l)+t\sum_l\beta_{ij}^l(v_l+tJv_l)=$$
$$\sum_l(\alpha_{ij}^l+t\beta_{ij}^l)(v_l+tJv_l).$$

\begin{nota}\rm
Given a real-analytic manifold $M$, let $M^\C$ denote the
germ along $M$ of a complexification of $M$, i.e. a
complex-analytic manifold $Z$ such that $M$ is
embedded in $Z$ as a completely real submanifold. The germ
$M^\C$ is defined uniquely up to (a germ of) a
biholomorphic map (see ~\cite{bw}). If $\phi$ is a
real-analytic function on $M$, $\phi^c$ will stand for the
unique germ along $M$ of a complex-analytic function on
$M^\C$ such that $\phi^c|_M=\phi$.
\end{nota}

\begin{thm}
Let $M^{k(n+1)}$ be a manifold of dimension $m=k(n+1)$ and
let $\{ w(t)\}_{t\in\R\P^1}$ be an integrable generalized
Veronese curve of distributions of codimension $k$ on $M$.
Then for any point $x\in M$ there exist a coordinate map
$(U,\phi), M\supset U\ni x, \phi=(\phi_1,\ldots,\phi_m)$
and a germ along $U$ of an integrable distribution
$F\subset TU^\C$ ($T$ stands for the real tangent bundle)
 such that for any $t\in\R\P^1$ :
\begin{enumerate}
\item[(1)] the distribution $(\Id+tJ)F\subset
TU^\C$ is integrable (we assume that the value $t=\infty$
corresponds to $JF$);
\item[(2)] $\ra((\Id+tJ)F)=kn=\ra w(t)$;
\item[(3)] the distribution $(\Id+tJ)F$ is projectable on $U$
along the germ of the foliation ${\cal
Y}=\{\Re\phi^c_1=\mbox{const},\ldots\Re\phi^c_m=\mbox{const}\}$;
\item[(4)] the projection of $(\Id-tJ)F$ coincides with
$w(t)|_U$.
\end{enumerate}
\end{thm}

\noindent{\sc Proof} Let $a_1,\ldots, a_{n+1}\in\R$ be
different nonzero numbers. Then the foliations ${\cal
W}(a_1),\ldots,\lbr {\cal W}(a_{n+1})$ are in general
position and for any point one can find a neighbourhood $U$
and functions
$$
\begin{array}{rcr}
\phi_1=\psi_1^1 &
,\ldots, & \phi_{n+1}=\psi_{n+1}^1 \\
\phi_{(n+1)+1}=\psi_1^2 & ,\ldots, &
\phi_{2n+2}=\psi_{n+1}^2 \\
 & \ldots & \\
\phi_{(k-1)(n+1)+1}=\psi_1^k & ,\ldots, &
\phi_m=\psi_{n+1}^k
\end{array}
$$
on $U$ such that ${\cal
W}(a_j)=\{\psi_j^1=\mbox{const},\ldots,\psi_j^k=\mbox{const}\},
j\in\overline{1,n+1}$.

We define a new family of distributions
$F(t)\subset TU^\C, t\in\R\P^1$, of codimension $k$
 by
$$
\G
(F(t))^\bot=\langle (\Id-tJ^*)\pi^*\g^1(t),\ldots
,(\Id-tJ^*)\pi^*\g^k(t)\rangle .
$$
Here $\langle\cdot\rangle$ stands for the linear span,
$J^*:T^*U^\C\rightarrow T^*U^\C$ is the adjoint operator to
the complex structure $J:TU^\C\rightarrow TU^\C$,
$\pi:U^\C\rightarrow U$ is the projection along the
foliation ${\cal Y}$ defined in (3), and
$\g^1(t),\ldots,\g^k(t)$ are the annihilating $w(t)$
1-forms (see Definition \ref{d1}).

Now, let us define the distribution $F$ as
$$F=\bigcap_{t\in\R\P^1}F(t)
$$
or
$$
\G
F^\bot=\langle (\Id-tJ^*)\pi^*\g^1(t),\ldots,
(\Id-tJ^*)\pi^*\g^k(t)|t\in\R\P^1\rangle .
$$

Notice that the 1-forms
$\pi^*\g_0^i,\pi^*\g_1^i-J^*\pi^*\g_0^i,\ldots,
\pi^*\g_n^i-J^*\pi^*\g_{n-1}^i,-J^*\pi^*\g_n^i$
corresponding to different powers of $t$ in
$(\Id-tJ^*)\pi^*\g^i(t), i\in\overline{1,k}$, are linearly
independent. Therefore the standard properties of the
Veronese curve (of degree $n+2$) imply that
$$
\G
F^\bot=\langle (\Id-tJ^*)\pi^*\g^1(t),\ldots,(\Id-tJ^*)\pi^*\g^k(t)|
t\in\{0,a_1,\ldots,a_{n+1}\}\rangle
$$
or
$$
F=\bigcap_{j\in\overline{0,n+1}}F(a_j)
$$
(here we put $a_0=0$).

This allow us to  prove the integrability of $F$ by
showing the integrability of $F(a_0),\ldots,\lbr
F(a_{n+1})$.

Evidently,
$\g^i(a_j)=\sum_{s=1}^k\beta_{js}^id\psi_j^s,i\in\overline{1,k},
j\in\overline{1,n+1}$ for some functions $\beta_{js}^i$.
Similarly,
$\g^i(a_0)=\sum_{s=1}^k\beta_{0s}^id\psi_0^s,i\in\overline{1,k}$
for some functions $\beta_{0s}^i,\psi_0^s$. Thus
\begin{eqnarray*}
(\Id-a_jJ^*)\pi^*\g^i(a_j)=
(\Id-a_jJ^*)\sum_{s=1}^k\pi^*\beta_{js}^id\pi^*\psi_j^s=\\
\sum_{s=1}^k\pi^*\beta_{js}^i(d(\Re(\psi_j^s)^c)+a_jd(\Im\psi_j^s)^c))=\\
\sum_{s=1}^k\pi^*\beta_{js}^id(\Re(\psi_j^s)^c+a_j(\Im\psi_j^s)^c),
i\in\overline{1,k}, j\in\overline{1,n+1}\\
\end{eqnarray*}
and
$$
(\Id-a_0J^*)\pi^*\g^i(a_0)=\pi^*\g^i(a_0)
=\sum_{s=1}^k\pi^*\beta_{0s}^id\pi^*\psi_0^s,
i\in\overline{1,k}
$$
(here we used the obvious facts that
$\pi^*\phi_j=\Re\phi_j^c, j\in\overline{1,m}$, and that
$J^*d(\Re\phi_j^c)=-d(\Im\phi_j^c)$). Now it is easy to
check the Frobenius integrability conditions using the
nondegeneracy of the matrix $\| \beta_{js}^i\|_{i,s}$ for
any $j\in\overline{0,n+1}$. So the distributions $F(a_j)$
are indeed integrable.

To prove (1) we choose another set of generators for
$\G F^\bot$ as follows
$$\G
F^\bot=\langle (\Id-tJ^*)\pi^*\g^1(t),\ldots,(\Id-tJ^*)\pi^*\g^k(t)|
t\in\{a_1,\ldots,a_{n+2}\}\rangle ,
$$
where $a_{n+2}:=\infty$, and notice that
$$
\G
J^*F^\bot=\langle (J^*+t\Id)\pi^*\g^1(t),\ldots,(J^*+t\Id)\pi^*\g^k(t)|
t\in\{a_1,\ldots,a_{n+2}\}\rangle .
$$
Now the integrability for
$JF=((J^*)^{-1}F^\bot)^\bot=(J^*F^\bot)^\bot$ can be proved
by the same considerations as for $F$ and the integrability
of $(\Id+tJ)F$ follows from Lemma.

In order to prove (2) we mention that $\ra F=kn$ by the
construction and that $\Id+tJ:TU^\C\rightarrow TU^\C$ is the
isomorphism for any $t\in\R\P^1$: the inverse operator is
given by the formula
$$
(\Id+tJ)^{-1}=\frac{1}{1+t^2}(\Id-tJ), t\not=\infty,
$$
$$
J^{-1}=-J.
$$

Now we are able to prove (3). We need to show that the
distribution $(\Id-tJ)F+T{\cal Y}\subset TU^\C$ is
integrable for any $t\in\R\P^1$. We fix $t=t_0\in\R$,
choose $b_0=t_0,b_1,\ldots,b_{n+1}$ to be different real
numbers and calculate the annihilators:
\begin{eqnarray*}
((\Id-t_0J)F+T{\cal Y})^\bot  =
 ((\Id-t_0J)F)^\bot\cap(T{\cal Y})^\bot, \\
\G(T{\cal
Y})^\bot  =  \langle d\Re\phi_1^c,\ldots,d\Re\phi_m^c\rangle , \\
\G((\Id-t_0J)F)^\bot  =  \G((\Id-t_0J^*)^{-1}F^\bot)
=\\
(\Id+t_0J^*)\langle(\Id-tJ^*)\pi^*\g^1(t),\ldots,(\Id-tJ^*)\pi^*\g^k(t)|
t\in\{b_0,\ldots,b_{n+1}\}\rangle   =  \\
\langle(1+t_0^2)\pi^*\g^1(t_0),\ldots,(1+t_0^2)\pi^*\g^k(t_0), \\
(1+t_0b_j)\pi^*\g^1(b_j)+(t_0-b_j)J^*\pi^*\g^1(b_j),\ldots,\\
(1+t_0b_j)\pi^*\g^k(b_j)+(t_0-b_j)J^*\pi^*\g^k(b_j)|j\in\overline{1,n+1}\rangle .
\end{eqnarray*}
It is easy to see that the collection of 1-forms
$$
\{\pi^*\g^i(b_j),J^*\pi^*\g^i(b_j)\}_{1\le i\le k, 1\le
j\le n+1}
$$
is a coframe on $U^\C$, hence the
1-forms $(t_0-b_j)J^*\pi^*\g^i(b_j),1\le i\le k, 1\le j\le
n+1$, can't be linearly expressed by the 1-forms
$d\Re\phi_1^c,\ldots,d\Re\phi_m^c$ which are combinations
of $\pi^*\g^i(b_j)$. So, finally,
$$
((\Id-t_0J)F+T{\cal
Y})^\bot=\langle \pi^*\g^1(t_0),\ldots,\pi^*\g^k(t_0)\rangle
$$
and the distribution
$$
(\Id-t_0J)F+T{\cal Y}=\pi^*w(t_0)
$$
is integrable. In the same manner one can show that
$$
JF+T{\cal Y}=\pi^*w(\infty).
$$

Simultaneously, the last two equations prove item (4).

\section{Application to integrability}

\begin{thm}
Let $M^{k(n+1)}$ be a manifold and let
$\{w(t)\}_{t\in\R\P^1}$ be a generalized Veronese curve of
distributions of codimension $k$ on $M$. Then in order that
$\{w(t)\}_{t\in\R\P^1}$  is a generalized Veronese web it
is sufficient that it is integrable at $n+3$ different
points $a_0,\ldots,a_{n+2}\in\R\P^1$.
\end{thm}

\noindent{\sc Proof} A careful analysis of the proof of
Theorem 1 shows that this proof uses only integrability
of $\{w(t)\}$ at $0,\infty$ and arbitrary
different nonzero finite points $a_1,\ldots,a_{n+1}$
for the construction of the
distribution $F\subset TU^\C$ such that $(\Id+tJ)F$ is
integrable for any $t\in\R\P^1$ and is projectable onto
$w(t)$. Thus it remains to map $0$ and $\infty$ to $a_0$
and $a_{n+2}$ respectively by an appropriate automorphism
of $\R\P^1$.

\end{document}